\documentclass[10pt,a4paper]{article}
\usepackage{graphicx}
\usepackage{theorem}
\usepackage{amsfonts}
\usepackage{hyperref}

\newcommand{\eps}{\varepsilon}
\newcommand{\disp}{\displaystyle}
\def\IN{\mathbb{N}}
\def\IZ{\mathbb{Z}}
\def\IR{\mathbb{R}}
\newcommand{\CC}{{\cal C}}

\newcommand{\subsetnoteq}{\subseteq_{\mbox{\rm \hspace{-0.5em}{\scriptsize /}
}}\hspace{-0.3em}}
\newcommand{\supsetnoteq}{\supseteq_{\mbox{\rm \hspace{-0.65em}{\scriptsize /}
}}\hspace{-0.15em}}

\newcommand{\carre}{\raise0.3em\hbox{\fbox{\rule{0ex}{0.3ex}%
\rule{0.3ex}{0ex}}}}

\newtheorem{theo}{Theorem}[section]
\newtheorem{prop}[theo]{Proposition}
\newtheorem{lem}[theo]{Lemma}
\newtheorem{cor}[theo]{Corollary}
\theorembodyfont{\rm}
\newtheorem{defi}[theo]{Definition}
\newtheorem{ex}[theo]{Example}
\newtheorem{rem}[theo]{Remark}

\newenvironment{demo}{\par\medskip\noindent 
{\it Proof}.}{\hfill\carre\par\medskip}

\newcommand{\lignelarge}{\lower1.5ex\hbox{\rule{0ex}{4ex}}}

\begin{document}

\title{On the Vere-Jones classification and existence of maximal measures 
for countable topological Markov chains
\footnotetext{Pacific J. Math., \textbf{209}, No. 2, 365--380, 2003.}
}
\author{Sylvie Ruette}
\date{}

\maketitle

\begin{abstract}
We consider topological Markov chains (also called Markov shifts) on
countable graphs.
We show that a transient graph can be extended to a recurrent graph of
equal entropy which is either positive recurrent of null recurrent, and we
give an example of each type.
We extend the notion of {\em local entropy} to topological 
Markov chains and prove that a
transitive Markov chain admits a measure of maximal entropy (or {\em maximal
measure}) whenever its local entropy is less than its (global) entropy.
\end{abstract}

\section*{Introduction}
In this article we are interested in connected oriented graphs and 
topological Markov
chains. All the graphs we consider have a countable set of 
vertices. If $G$ is an oriented graph, let $\Gamma_G$ be the set of
two-sided infinite sequences of vertices that form a path
in $G$ and let $\sigma$ denote the shift 
transformation. The {\em Markov chain} associated to $G$ is the (non
compact) dynamical system $(\Gamma_G,\sigma)$. The entropy $h(G)$ of the
Markov chain  $\Gamma_G$ was defined by Gurevich; it can be 
computed by several ways and satisfies the Variational Principle 
\cite{Gur1, Gur2}.

In \cite{Ver1} Vere-Jones classifies connected oriented graphs as transient,
null recurrent or positive recurrent according to the properties of the
series associated with the number of loops, by analogy with probabilistic
Markov chains. To a certain extent, positive
recurrent graphs resemble finite graphs. In \cite{Gur1} Gurevich 
shows that a Markov chain on a connected graph admits a measure of
maximal entropy (also called {\em maximal measure}) if and only if the
graph is positive recurrent. In this case, this measure is unique and it is
an ergodic Markov measure.

In \cite{Sal2,Sal3} Salama gives a geometric approach to the Vere-Jones 
classification. The fact that a graph can (or cannot)
be ``extended'' or ``contracted'' without changing its entropy is
closely related to its class. In particular a graph with no proper
subgraph of equal entropy is positive recurrent. The converse is not true
\cite{Sal3} (see also
\cite{uFie} for an example of a positive recurrent graph with a finite
valency at every vertex that has no proper subgraph of equal entropy). 
This result shows that the
positive recurrent class splits into two subclasses: a graph is called
{\em strongly positive recurrent} if it has no proper subgraph of equal 
entropy; it is equivalent to a combinatorial condition (a finite
connected graph is always strongly positive recurrent).
In \cite{Sal2, Sal3} Salama also states that a graph is transient if and only
if it can be extended to a bigger transient graph of equal entropy.
We show that any transient graph $G$ is contained in 
a recurrent graph of equal entropy, which is positive or null recurrent
depending on the properties of $G$. We illustrate the two possibilities
-- a transient graph with a positive or null recurrent extension -- 
by an example. 

The result of Gurevich entirely solves the question of existence of a
maximal measure in term of graph classification. Nevertheless it is not so
easy to prove that a graph is positive recurrent and one may wish to have 
more efficient criteria. In \cite{GZ2} Gurevich and Zargaryan give a
sufficient condition for existence of a maximal measure; 
it is formulated in terms of exponential
growth of the number of paths inside and outside a finite subgraph.
We give a new sufficient criterion based on local entropy.

Why consider local entropy? For a compact dynamical system, it
is known that a null local entropy implies the existence of a maximal
measure (\cite{New}, see also \cite{Bow4} for a similar but  
different result). This result may be strengthened in some cases: it is
conjectured that, if $f$ is a  map of the interval which is $C^r$, $r>1$,
and satisfies $h_{top}(f)>h_{loc}(f)$, 
then there exists a maximal measure \cite{Buz}.
Our initial motivation comes from the conjecture above because smooth interval
maps and Markov chains are closely related.
If $f\colon [0,1]\to [0,1]$ is $C^{1+\alpha}$ (i.e. $f$ is $C^1$ and
$f'$ is $\alpha$-H\"older with $\alpha>0$) with $h_{top}(f)>0$ then an oriented 
graph $G$ can be associated to $f$, $G$ is connected if $f$ is transitive,
  and there is a bijection between the maximal measures of $f$ and those 
of $\Gamma_G$ \cite{Buz, Buz4}. 
We show that a Markov chain is strongly positive recurrent, 
thus admits a maximal measure, if its local entropy is
strictly less that its Gurevich entropy.
However this result does not apply directly to interval maps
since the ``isomorphism'' between $f$ and its Markov extension is 
not continuous so it may not preserve local entropy (which
depends on the distance).

The article is organized as follows. Section~\ref{sec:background}
contains definitions and basic properties on oriented graphs and Markov
chains. In Section~\ref{sec:classification}, after recalling the
definitions of transient, null recurrent and positive recurrent graphs
and some related properties, we show that any transient graph is contained
in a recurrent graph of equal entropy 
(Proposition~\ref{prop:transient-recurrent}) 
and we give an example of a transient graph
which extends to a positive recurrent (resp. null recurrent) graph. 
Section~\ref{sec:maximal-measure} is devoted to
the problem of existence of maximal measures: 
Theorem~\ref{theo:hloc-maximal-measure}
gives a sufficient condition for the existence of a maximal measure, based on
local entropy. 

\section{Background}\label{sec:background}
\subsection{Graphs and paths}

Let $G$ be an oriented graph with a countable set of vertices $V(G)$. If $u,v$
are two vertices, there is at most one arrow $u\to v$. A {\em path} of 
length $n$ is a sequence of vertices 
$(u_0,\cdots,u_n)$ such that $u_i\to u_{i+1}$
in $G$ for $0\leq i<n$. This path is called a {\em loop} if $u_0=u_n$.
We say that the graph $G$ is {\em connected} if for all vertices $u,v$ 
there exists a path from $u$ to $v$; in the literature, such a graph
is also called {\em strongly connected}.

If $H$ is a subgraph of $G$, we write $H\subset G$; if in addition $H\not=G$,
we write $H\subsetnoteq G$ and say that $H$ is a {\em proper subgraph}.
If $W$ is a subset of $V(G)$, the set $V(G)\setminus W$ is denoted by 
$\overline{W}$. We also denote by $W$ the subgraph of $G$ whose vertices 
are $W$ and whose edges are all edges of $G$ between two vertices in $W$. 

\medskip
Let $u,v$ be two vertices. We define the following quantities.
\begin{itemize}
\item $p_{uv}^G(n)$ is the number of paths $(u_0,\cdots,u_n)$ 
such that $u_0=u$ and $u_n=v$; $R_{uv}(G)$ is the radius of convergence
of the series $\sum p_{uv}^G(n)z^n$.
\item $f_{uv}^G(n)$ is the number of paths $(u_0,\cdots,u_n)$ 
such that $u_0=u$, $u_n=v$ and $u_i\not=v$ for $0<i<n$; 
$L_{uv}(G)$ is the radius of convergence
of the series $\sum f_{uv}^G(n)z^n$.
\end{itemize}

\begin{prop}[Vere-Jones \cite{Ver1}]
Let $G$ be an  oriented graph. If $G$ is connected, 
$R_{uv}(G)$ does not depend on $u$ and $v$; it denoted by $R(G)$.
\end{prop}

If there is no confusion, $R(G)$ and $L_{uv}(G)$ will be written $R$ and 
$L_{uv}$. For a graph $G'$ these two radii will be written $R'$ and $L'_{uv}$.

\subsection{Markov chains} \label{subsec:markov-chain}
Let $G$ be an oriented graph.
$\Gamma_G$ is the set of two-sided infinite paths in $G$, that is,
$$
\Gamma_G=\{(v_n)_{n\in\IZ} \mid  \forall n\in \IZ, v_n\to v_{n+1} 
\mbox{ in } G \}\subset (V(G))^{\IZ}.
$$
$\sigma$ is the shift on $\Gamma_G$. The {\em (topological) Markov chain} 
on the graph $G$ is the system $(\Gamma_G,\sigma)$.

The set $V(G)$ is endowed with the discrete topology and $\Gamma_G$
is endowed with the induced topology of $(V(G))^{\IZ}$. 
The space $\Gamma_G$ is not compact unless $G$ is finite.
A compatible distance on $\Gamma_G$ is given by $d$,
defined as follows: $V(G)$ is identified with $\IN$ and the distance 
$D$ on $V(G)$ is given by $
D(n,m)=\left|\frac{1}{2^n}-\frac{1}{2^m}\right|$.
If $\bar u=(u_n)_{n\in\IZ}$ and $\bar v=(v_n)_{n\in\IZ}$ are two
elements of $\Gamma_G$,
$$
d(\bar u,\bar v)=\sum_{n\in\IZ}
\frac{D(u_n,v_n)}{2^{|n|}}\leq 3.
$$

The Markov chain $(\Gamma_G,\sigma)$ is 
{\em transitive} if for any non empty open
sets $A,B\subset \Gamma_G$ there exists $n>0$ such that $\sigma^n(A)\cap
B\not=\emptyset$. Equivalently, $\Gamma_G$ is transitive if and only if
the graph $G$ is connected. In the sequel we will be interested in
connected graphs only.

\subsection{Entropy}\label{subsec:def-entropy}

If $G$ is a finite graph, $\Gamma_G$ is compact and
the topological entropy $h_{top}(\Gamma_G,\sigma)$ is well defined
(see e.g. \cite{DGS} for the definition of the topological entropy).
If $G$ is a countable graph, the {\em Gurevich entropy} \cite{Gur1} 
of $G$ is given by
$$
h(G)=\sup\{h_{top}(\Gamma_H,\sigma)\mid H\subset G, H \mbox{ finite}\}.
$$

This entropy can also be computed in a combinatorial way, as the exponential
growth of the number of paths with fixed endpoints \cite{Gur2}.

\begin{prop}[Gurevich]
Let $G$ be a connected oriented graph. Then for any vertices $u,v$
$$
h(G)=\lim_{n\to+\infty}\frac{1}{n}\log p_{uv}^G(n)=-\log R(G).
$$
\end{prop}

Another way to compute the entropy is to compactify the space $\Gamma_G$ and
then use the definition of topological entropy for compact metric spaces.
If $G$ is an oriented graph, denote
the one-point compactification of $V(G)$ by $V(G)\cup\{\infty\}$ and
define $\overline{\Gamma}_G$ as the closure of $\Gamma_G$ in 
$(V(G)\cup\{\infty\})^{\IZ}$. The distance $d$ naturally
extends to $\overline{\Gamma}_G$.
In \cite{Gur1} Gurevich shows that this gives the
same entropy; this
means that there is only very little dynamics added in this compactification.
Moreover, the Variational Principle is still valid for Markov chains
\cite{Gur1}.

\begin{theo}[Gurevich]
Let $G$ be an oriented graph. Then
$$
h(G)=h_{top}(\overline{\Gamma}_G,\sigma)=\sup\{h_{\mu}(\Gamma_G)\mid \mu
\ \sigma\mbox{-invariant probability measure}\}.
$$
\end{theo}

\section{On the classification of connected graphs}\label{sec:classification}
\subsection{Transient, null recurrent, positive recurrent graphs}
In \cite{Ver1} Vere-Jones gives a classification of connected graphs
as transient, null recurrent or positive recurrent. The definitions
are given in Table~\ref{tab:classification} (lines 1 and 2) as well as
properties of the series $\sum p_{uv}^G(n)z^n$ which give an alternative
definition. 

\begin{table}[ht] 
\begin{center}
\begin{tabular}{l|c|c|c}
                     & transient   & null      & positive  \\
                     &             & recurrent & recurrent \\
\hline
\lignelarge
$\disp\sum_{n>0} f^G_{uu}(n)R^n$  & $<1$        & $1$       & $1$     \\
\hline
\lignelarge
$\disp\sum_{n>0} nf^G_{uu}(n)R^n$ &$\leq+\infty$&$+\infty$  &$<+\infty$ \\
\hline
\lignelarge
$\disp\sum_{n\geq 0} p^G_{uv}(n)R^n$  &$<+\infty$   &$+\infty$  &$+\infty$ \\
\hline
\lignelarge
$\disp\lim_{n\to+\infty} p^G_{uv}(n)R^n$  & $0$         & $0$       &$\lambda_{uv}>0$\\
\hline
\lignelarge
                     & $R=L_{uu}$       & $R=L_{uu}$     & $R\leq L_{uu}$
\end{tabular}
\end{center}
\caption{properties of the series associated to a transient, null 
recurrent or positive recurrent graph $G$; these
properties do not depend on the vertices $u,v$ ($G$ is connected). \label{tab:classification}}
\end{table}

In \cite{Sal2,Sal3} Salama studies the links between the classification
and the possibility to extend or contract a graph
without changing its entropy. It follows that a connected graph is transient
if and only if it is strictly included in a connected graph of equal
entropy, and that a graph with no proper subgraph of equal entropy 
is positive recurrent.

\begin{rem}
In \cite{Sal2} Salama claims that $L_{uu}$ is independent of $u$, which is
not true;  in \cite{Sal3} he
uses the quantity $L=\inf_u L_{uu}$ and he states that if $R=L$ then 
$R=L_{uu}$ for all vertices $u$, which is wrong too (see Proposition~3.2 in
\cite{GurS}). It follows that in \cite{Sal2,Sal3} the statement ``R=L'' 
must be interpreted
either as ``$R=L_{uu}$ for some $u$'' or ``$R=L_{uu}$ for all $u$''
depending on the context. 
This encouraged us to give the proofs of Salama's results in this article.
\end{rem}

In \cite{Sal3}
Salama shows that a transient or null recurrent graph satisfies $R=L_{uu}$
for all vertices $u$; we give the unpublished proof due to
U. Fiebig \cite{uFie}.

\begin{prop}[Salama] \label{prop:transient-R=L}
Let $G$ be a connected oriented graph. If $G$ is transient or null
recurrent then $R=L_{uu}$ for all vertices $u$. Equivalently,
if there exists a vertex $u$ such that $R<L_{uu}$ then $G$ is positive 
recurrent.
\end{prop}

\begin{demo}
For a connected oriented graph, it is obvious that $R\leq L_{uu}$ for all
$u$, thus the two claims of the Proposition are equivalent. We
prove the second one.

Let $u$ be a vertex of $G$ such that $R<L_{uu}$. Let
$F(x)=\sum_{n\geq 1}f_{uu}^G(n)x^n$ for all $x\geq 0$. If we break a loop 
based in $u$ into first return loops, we get the following formula:
\begin{equation}\label{eq:P(x)-F(x)}
\sum_{n\geq 0}p_{uu}^G(n)x^n=\sum_{k\geq 0}(F(x))^k.
\end{equation}
Suppose that $G$ is transient,
that is, $F(R)<1$. The map $F$ is analytic on $[0,L_{uu})$ and $R<L_{uu}$
thus there exists $R<x<L_{uu}$ such that $F(x)<1$. According 
to Equation~(\ref{eq:P(x)-F(x)})
one gets that $\sum_{n\geq 0}p_{uu}^G(n)x^n<+\infty$, which contradicts
the definition of $R$. Therefore $G$ is recurrent. Moreover $R<L_{uu}$
by assumption, thus $\sum_{n\geq 1} nf_{uu}^G(n)R^n<+\infty$, which
implies that $G$ is positive recurrent.
\end{demo}

\begin{defi}
A connected oriented graph is called {\em strongly positive recurrent} if
$R<L_{uu}$ for all vertices $u$.
\end{defi}

\begin{lem}\label{lem:recurrent-R}
Let $G$ be a connected oriented graph and $u$ a vertex.
\begin{enumerate}
\item $R<L_{uu}$ if and only if $\,\sum_{n\geq 1}f_{uu}^G(n)L_{uu}^n>1$.
\item If $G$ is recurrent then $R$ is the unique positive number $x$ such that
$\sum_{n\geq 1}f_{uu}^G(n)x^n=1$.
\end{enumerate}
\end{lem}

\begin{demo} Use the fact that 
$F(x)=\sum_{n\geq 1}f_{uu}^G(n)x^n$ is increasing.
\end{demo}

The following result deals with transient graphs \cite{Sal2}.

\begin{theo}[Salama]
Let $G$ be a connected oriented graph of finite positive entropy.
Then $G$ is transient if and only if there exists a connected oriented
graph $G'\supsetnoteq G$ such that $h(G')=h(G)$. If $G$ is transient then
$G'$ can be chosen transient.
\end{theo}

\begin{demo}
The assumption on the entropy implies that $0<R<1$.
Suppose first that there exists a connected graph $G'\supsetnoteq G$ such
that $h(G')=h(G)$, that is, $R'=R$. Fix a vertex $u$ in $G$. 
The graph $G$ is a proper subgraph of $G'$ thus there exists $n$ such that 
$f_{uu}^G(n)<f_{uu}^{G'}(n)$, which implies that
$$
\sum_{n\geq 1}f_{uu}^G(n)R^n<\sum f_{uu}^{G'}(n){R'}^n\leq 1.
$$
Therefore $G$ is transient.

Now suppose that $G$ is transient and fix a vertex $u$ in $G$. One has
$\sum_{n\geq 1} f_{uu}^G(n)R^n<1$. Let $k\geq 2$ be an integer such that
$$
\sum_{n\geq 1} f_{uu}^G(n)R^n+R^k<1.
$$
Define the graph $G'$ by adding a loop of length $k$ based at the vertex 
$u$; one has $R'\leq R$ and
\begin{equation}\label{eq:G'-G-transient}
\sum_{n\geq 1}f_{uu}^{G'}(n){R'}^n\leq 
\sum_{n\geq 1}f_{uu}^{G'}(n){R}^n=
\sum_{n\geq 1}f_{uu}^{G}(n){R}^n+R^k<1.
\end{equation}
Equation~(\ref{eq:G'-G-transient}) implies that $R\leq L_{uu}'$ and also
that the graph $G'$ is transient, so $R'=L_{uu}'$ by
Proposition~\ref{prop:transient-R=L}. 
Then one has $L_{uu}'=R'\leq R\leq L_{uu}'$ thus $R=R'$.
\end{demo}

In \cite{Sal3} Salama proves that if $R=L_{uu}$ for all vertices $u$ then
there exists a proper subgraph of equal entropy. We show that the
same conclusion holds if one supposes that $R=L_{uu}$ for some $u$. The
proof below is a variant of the one of Salama. The converse is also true,
as shown by U. Fiebig \cite{uFie}. 

\begin{prop}\label{prop:R<=L}
Let $G$ be a connected oriented graph of positive entropy.
\begin{enumerate}
\item If there is a vertex $u$ such that $R=L_{uu}$ then there exists
a connected subgraph $G'\subsetnoteq G$ such that $h(G')=h(G)$.
\item If there is a vertex $u$ such that $R<L_{uu}$ then for all proper
subgraphs $G'$ one has $h(G')<h(G)$.
\end{enumerate}
\end{prop}

\begin{demo}
\par\noindent i) Suppose that $R=L_{uu}$. If 
$u_0=u$ is followed by a unique vertex, let $u_1$ be this vertex. If
$u_1$ is followed by a unique vertex, let $u_2$ be this vertex, and so on.
If this leads to define $u_n$ for all $n$ then $h(G)=0$, which is not
allowed.

Let $u_k$ be the last built vertex; there exist two distinct vertices
$v,v'$ such that $u_k\to v$ and $u_k\to v'$. Let $G_1'$ be the graph
$G$ deprived of the arrow $u_k\to v$ and $G_2'$ the graph $G$ deprived
of all the arrows $u_k\to w$, $w\not=v$. Call $G_i$ the connected component
of $G_i'$ that contains $u$ ($i=1,2$); obviously $G_i\subsetnoteq G$. 
For all $n\geq 1$ one has
$$
f_{uu}^G(n)=f_{u_ku}^G(n-k)=f_{u_ku}^{G_1}(n-k)+f_{u_ku}^{G_2}(n-k),
$$
thus there exists $i\in\{1,2\}$ such that $L_{uu}=L_{u_ku}^{G_i}$. One has
$$
R\leq R(G_i)\leq L_{u_ku}(G_i)=L_{uu}=R,
$$
thus $R=R(G_i)$, that is, $h(G)=h(G_i)$.

\medskip\noindent ii)
Suppose that $R<L_{uu}$ and consider $G'\subsetnoteq G$. Suppose first
that $u$ is a vertex of $G'$. The graph $G$ is positive recurrent by
Proposition~\ref{prop:transient-R=L} so $\sum_{n\geq 1}f_{uu}^G(n)R^n=1$.
Since $G'\subsetnoteq G$ there exists $n$ such that 
$f_{uu}^{G'}(n)<f_{uu}^G(n)$, thus
\begin{equation}\label{eq:G'-G-recurrent}
\sum_{n\geq 1} f_{uu}^{G'}R^n<1.
\end{equation}
Moreover $L_{uu}'\geq L_{uu}$. If $G'$ is transient then $R'=L_{uu}'$
(Proposition~\ref{prop:transient-R=L}) thus $R'\geq L_{uu}>R$. If
$G'$ is recurrent then $\sum_{n\geq 1}f_{uu}^{G'}{R'}^n=1$ thus $R'>R$
because of Equation~(\ref{eq:G'-G-recurrent}). In both cases $R'>R$,
that is, $h(G')<h(G)$.

Suppose now that $u$ is not a vertex of $G'$ and fix a vertex $v$ in $G'$.
Let $(u_0,\ldots,u_p)$ a path (in $G$) of minimal length between $u=u_0$
and $v=u_p$, and let $(v_0,\ldots, v_q)$ be a path of minimal length
between $v=v_0$ and $u=v_q$. 

If $(w_0=v,w_1,\ldots, w_n=v)$ is a loop in $G'$ then
$$
(u_0=u,u_1,\ldots,u_p=w_0,w_1,\ldots,w_n=v_0,v_1,\ldots,v_q=u)
$$
is a first return loop based in $u$ in the graph $G$. For all $n\geq 0$
we get that $p_{vv}^{G'}(n)\leq f_{uu}^G(n+p+q)$, thus
$R'\geq L_{uu}>R$, that is, $h(G')<h(G)$.
\end{demo}

The following result gives a characterization of strongly positive
recurrent graphs. It is a straightforward corollary of 
Proposition~\ref{prop:R<=L} (see also \cite{uFie}).

\begin{theo}\label{theo:strongly-positive-recurrent}
Let $G$ be a connected oriented graph of positive entropy. The following
properties are equivalent:
\begin{enumerate}
\item for all $u$ one has $R<L_{uu}$ (that is, $G$ is strongly positive 
recurrent),
\item there exists $u$ such that $R<L_{uu}$,
\item $G$ has no proper subgraph of equal entropy.
\end{enumerate}
\end{theo}

\subsection{Recurrent extensions of equal entropy of transient graphs}
\label{subsec:exR=L}

We show that any transient graph $G$ can be extended to a
recurrent graph without changing the entropy by adding a (possibly infinite)
number of loops. If the series $\sum n f_{vv}^G(n)R^n$ is finite then 
the obtained recurrent graph is positive recurrent (but not strongly 
positive recurrent), otherwise it is null recurrent.

\begin{prop}\label{prop:transient-recurrent}
Let $G$ be a transient graph of finite positive entropy. Then there exists a
recurrent graph $G'\supset G$ such that $h(G)=h(G')$. Moreover $G'$ can be
chosen to be positive recurrent if $\,\sum_{n>0} nf_{uu}^G(n)R^n<+\infty$
for some vertex $u$ of $G$, and $G'$ is necessarily null recurrent otherwise.
\end{prop}

\begin{demo}
The entropy of $G$ is finite and positive thus $0<R<1$ 
and there exists an integer $p$ such that $\frac{1}{2}\leq pR<1$.
Define $\alpha=pR$.
Let $u$ be a vertex of $G$ and define 
$D=1-\sum_{n\geq 1} f_{uu}^G(n)R^n$; one has $0<D<1$. Moreover
$$
\sum_{n\geq 1} \alpha^n\geq \sum_{n\geq 1} \frac{1}{2^n}=1,
$$
thus
\begin{equation}\label{eq:series}
\sum_{n\geq k+1}\alpha^n=\alpha^k\sum_{n\geq 1}\alpha^n\geq \alpha^k.
\end{equation}

We build a sequence of integers $(n_i)_{i\in I}$ such that
$2\sum_{i\in I} \alpha^{n_i}=D$. For this, we define
inductively a strictly increasing (finite or infinite)
sequence of integers $(n_i)_{i\in I}$ such that for all $k\in I$ 
$$
\sum_{i=0}^k \alpha^{n_i} \leq \frac{D}{2} < 
\sum_{i=0}^k \alpha^{n_i} + \sum_{n>n_k} \alpha^n.
$$

\noindent
-- Let $n_0$ be the greatest integer $n\geq 2$ such that
$\sum_{k\geq n} \alpha^k >\frac{D}{2}$. By choice of $n_0$ one has
$\sum_{n\geq n_0+1} \alpha^n \leq \frac{D}{2}$,
thus $\alpha^{n_0}\leq\frac{D}{2}$ by Equation~(\ref{eq:series}). 
This is the required property at rank $0$.

\noindent
-- Suppose that $n_0,\cdots,n_k$ are already defined. 
If $\sum_{i=0}^k \alpha^{n_i}=
\frac{D}{2}$ then $I=\{0,\cdots,k\}$ and we stop the construction.
Otherwise let $n_{k+1}$ be the
greatest integer $n >n_k$ such that 
$$
\sum_{i=0}^k \alpha^{n_i}+\sum_{j\geq n} \alpha^j > \frac{D}{2}.
$$
By choice of $n_{k+1}$ and Equation~(\ref{eq:series}), one has
$$ 
\alpha^{n_{k+1}}\leq \sum_{j\geq n_{k+1}+1} \alpha^j\leq 
\frac{D}{2}-\sum_{i=0}^k \alpha^{n_i}.
$$
This is the required property at rank
$k+1$.

\medskip
Define a new graph $G'\supset G$ by adding $2p^{n_i}$ loops of length
$n_i$ based at the vertex $u$. Obviously one has $R'\leq R$, and
$\sum_{i\in I}(pR)^{n_i}=\frac{D}{2}$ by construction.
Therefore
\begin{equation}\label{eq:sum=1}
\sum_{n\geq 1} f_{uu}^{G'}(n)R^n=
\sum_{n\geq 1} f_{uu}^G(n)R^n+\sum_{i\in I}2(pR)^{n_i}=1.
\end{equation}
This implies that $R\leq L_{uu}'$. If $G'$ is transient then
$\sum_{n\geq 1} f_{uu}^{G'}(n){R'}^n<1$ and $R'=L_{uu}'$ by
Proposition~\ref{prop:transient-R=L}, thus $R\leq R'$ and
Equation~(\ref{eq:sum=1}) leads to a contradiction. Therefore $G'$ is
recurrent. By Lemma~\ref{lem:recurrent-R}(ii) one has $R'=R$, that is,
$h(G')=h(G)$.
In addition, 
$$
\sum_{n\geq 1} n f_{uu}^{G'}(n)R^n=
\sum_{n\geq 1} n f_{uu}^G(n)R^n+\sum_{i\in I} n_i \alpha^{n_i}
$$
and this quantity is finite if and only if $\sum n f_{uu}^G(n)R^n$ is
finite. In this case the graph $G'$ is positive recurrent.

\medskip
If $\sum n f_{uu}^G(n)R^n=+\infty$, let $H$ be a recurrent graph containing 
$G$ with $h(H)=h(G)$. Then $H$ is null recurrent because
$$
\sum_{n\geq 1} n f_{uu}^{H}(n)R^n\geq \sum_{n\geq 1} n f_{uu}^G(n)R^n=+\infty.
$$
\end{demo}

\begin{ex}\label{ex:R=L}
We build a positive (resp. null) recurrent graph $G$ such that
$\sum f_{uu}^G(n)L_{uu}^n=1$ and then we delete an arrow to obtain a
graph $G'\subset G$ which is transient and such that $h(G')=h(G)$. 
First we give a
description of $G$ depending on a sequence of integers $a(n)$ then we give
two different values to the sequence $a(n)$ so as to obtain a positive 
recurrent graph in one case and a null recurrent graph in the other case.

Let $u$ be a vertex and $a(n)$ a sequence of non negative integers 
for $n\geq 1$, with $a(1)=1$.
The graph $G$ is composed of $a(n)$ loops of length $n$ based at the vertex
$u$ for all $n\geq 1$ (see Figure~\ref{fig:G-G'}). More precisely, 
define the set of vertices of $G$ as 
$$
V=\{u\}\cup\bigcup_{n=1}^{+\infty}\{v_k^{n,i} \mid 1\leq i\leq a(n), 
1\leq k\leq n-1\},
$$
where the vertices $v_k^{n,i}$ above are distinct.
Let $v_0^{n,i}=v_n^{n,i}=u$ for $1\leq i\leq a(n)$. There is an arrow
$v_k^{n,i}\to v_{k+1}^{n,i}$ for $0\leq k\leq n-1, 1\leq i\leq a(n), n\geq 1$ 
and there is no other arrow in $G$. The graph $G$ is connected  
and $f_{uu}^G(n)=a(n)$ for $n\geq 1$. 

\begin{figure}[ht]
\centerline{\includegraphics{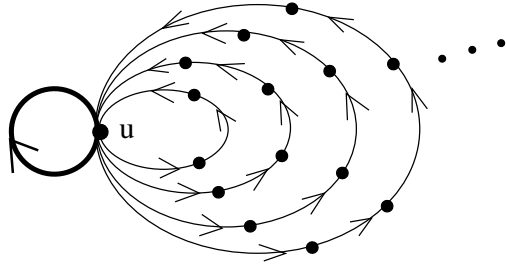}}
\caption{the graphs $G$ and $G'$; the bold loop (on the left)
is the only arrow that belongs to $G$ and not to $G'$, otherwise the
two graphs coincide. \label{fig:G-G'}}
\end{figure}

The sequence $(a(n))_{n\geq 2}$ is chosen such that it satisfies
\begin{equation}\label{eq:=1}
\sum_{n\geq 1} a(n)L^n=1,
\end{equation}
where $L=L_{uu}>0$ is the radius of convergence of the series $\sum a(n)z^n$.
If $G$ is transient then $R=L_{uu}$ by 
Proposition~\ref{prop:transient-R=L}, but 
Equation~(\ref{eq:=1}) contradicts the definition of transient. Thus
$G$ is recurrent. Moreover, $R=L$ by Lemma~\ref{lem:recurrent-R}(ii).

The graph $G'$ is obtained from $G$ by deleting the arrow $u\to u$.
Obviously one has $L_{uu}'=L$ and
$$
\sum_{n\geq 1} f_{uu}^{G'}(n)L^n=1-L<1.
$$
This implies that $G'$ is transient because $R'\leq L_{uu}'$. Moreover
$R'=L_{uu}'$ by Proposition~\ref{prop:transient-R=L}
thus $R'=R$, that is, $h(G')=h(G)$.

\medskip
Now we consider two different sequences $a(n)$.

\medskip\noindent {\bf 1)}
Let $a(n^2)=2^{n^2-n}$ for $n\geq 1$ and $a(n)=0$ otherwise.
Then $L=\frac{1}{2}$ and
$$
\sum_{n\geq 1}f_{uu}^G(n)L^n=\sum_{n\geq 1}2^{n^2-n}\frac{1}{2^{n^2}}=
\sum_{n\geq 1}\frac{1}{2^n}=1.
$$
Moreover
$$
\sum_{n\geq 1}nf_{uu}^G(n)L^n=\sum_{n\geq 1} \frac{n^2}{2^n}<+\infty,
$$
hence the graph $G$ is positive recurrent. 

\medskip\noindent {\bf 2)}
Let $a(1)=1$,  $a(2^n)=2^{2^n-n}$ for $n\geq 2$ 
and $a(n)=0$ otherwise. One can compute that $L=\frac{1}{2}$, and
$$
\sum_{n\geq 1}f_{uu}^G(n)L^n=
\frac{1}{2}+ \sum_{n\geq 2}2^{2^n-n}\frac{1}{2^{2^n}}
=\frac{1}{2}+\sum_{n\geq 2}\frac{1}{2^n}
=1.
$$
Moreover
$$
\sum_{n\geq 1}nf_{uu}^G(n)L^n
=\frac{1}{2}+\sum_{n\geq 2} 2^n\frac{1}{2^n}=+\infty
$$
hence the graph $G$ is null recurrent.
\end{ex}

\begin{rem}
Let $G$ be transient graph of finite entropy. Fix 
a vertex $u$ and choose an integer $k$ such that
$\sum_{n\geq k}R^n<1-\sum_{n\geq 1}f_{uu}^G(n)R^n$. 
For every integer $n\geq k$ let $m_n=\lfloor R^{-n}\rfloor$, add
$\lfloor R^{-(m_n-n)}\rfloor$ loops of length $m_n$ based at the vertex $u$
and call $G'$ the graph obtained in this way. It can be shown that
the graph $G'$ is transient, $h(G')=h(G)$ and $\sum_{n\geq 1}
nf_{uu}^{G'}R'^n=+\infty$. Then
Proposition~\ref{prop:transient-recurrent}  implies that every transient
graph is included in a null recurrent graph of equal entropy.
\end{rem}

\begin{rem}
In the more general setting of thermodynamic formalism for
countable Markov chains, Sarig puts to the fore
a subclass of positive recurrent potentials which he calls
strongly positive recurrent \cite{Sar4}; his
motivation is different, but the classifications agree.
If $G$ is a countable oriented graph, a {\em potential} is a
continuous map $\phi\colon \Gamma_G\to \IR$ and the {\em pressure}
$P(\phi)$ is the analogous of the Gurevich entropy, the paths being
weighted by $e^{\phi}$; a potential is either transient or null recurrent or
positive recurrent. Considering the null potential $\phi\equiv 0$, we 
retrieve the case of (non weighted) topological
Markov chains. In \cite{Sar4} Sarig introduces a quantity $\Delta_u[\phi]$; 
$\phi$ is transient (resp. recurrent) if
$\Delta_u[\phi]<0$ (resp. $\Delta_u[\phi]\geq 0$). The potential is called
{\em strongly positive recurrent} if $\Delta_u[\phi]>0$, which implies it is
positive recurrent. A strongly positive recurrent
potential $\phi$ is stable under perturbation, that is, any 
potential $\phi+t\psi$ close to $\phi$ is positive recurrent too. 
For the null potential,
$\Delta_u[0]=\log\left(\sum_{n\geq 1} f_{uu}(n)L^n\right)$, thus 
$\Delta_u[0]>0$ if and only if the graph is strongly positive recurrent
(Lemma~\ref{lem:recurrent-R} and 
Theorem~\ref{theo:strongly-positive-recurrent}).
In \cite{GurS} strongly positive recurrent potentials are called
{\em stable positive}.

Examples of (non null) potentials which are positive recurrent but not
strongly positive recurrent can be found in \cite{Sar4}; some of them resemble
much the Markov chains of Example~\ref{ex:R=L}, their graphs being 
composed of loops as in Figure~\ref{fig:G-G'}.
\end{rem}


\section{Existence of a maximal measure}\label{sec:maximal-measure}
\subsection{Positive recurrence and maximal measures}
A Markov chain on a finite graph always has a maximal measure \cite{Par},
but it is not the case for infinite graphs \cite{Gur1}. In \cite{Gur2} 
Gurevich gives a necessary and sufficient condition for the existence 
of such a measure.

\begin{theo}[Gurevich]\label{theo:Gurevic}
Let $G$ be a connected oriented graph of finite positive entropy. Then the
Markov chain $(\Gamma_G,\sigma)$ admits a maximal measure
if and only if the graph is positive recurrent. Moreover, such a measure
is unique if it exists, and it is an ergodic Markov measure.
\end{theo}

In \cite{GZ2} Gurevich and Zargaryan show that
if one can find a finite connected subgraph 
$H \subset G$ such that there are more paths inside than outside $H$
(in term of exponential growth), then the graph $G$ has a maximal
measure.
This condition is equivalent to strong positive recurrent as
it was shown by Gurevich and Savchenko in the more general setting of
weighted graphs~\cite{GurS}. 

\medskip
Let $G$ be a connected oriented graph, $W$ a subset of vertices 
and $u,v$ two vertices of $G$. Define
$t_{uv}^W(n)$ as the number of paths $(v_0,\cdots,v_n)$ such that $v_0=u,
v_n=v$ and $v_i\in W$ for all $0<i<n$, and put
$\disp \tau_{uv}^W=\limsup_{n\to+\infty} \frac{1}{n}\log t_{uv}^W(n)$.

\begin{theo}[Gurevich-Zargaryan]\label{theo:GZ-et-reciproque}
Let $G$ be a connected oriented graph of finite positive entropy.
If there exists a finite set of vertices $W$ such that $W$ is
connected and for all vertices $u,v$ in $W$, 
$\tau_{uv}^{\overline{W}}\leq h(W)$, then the graph
$G$ is strongly positive recurrent.
\end{theo}

For graphs that are not strongly positive recurrent 
the entropy is mainly
concentrated near infinity in the sense that it is supported by the infinite
paths that spend most of the time outside a finite subgraph
(Proposition~\ref{prop:entropy-infinity}). This result is obtained by applying
inductively the construction of Proposition~\ref{prop:R<=L}(i). As a
corollary, there exist ``almost maximal measures escaping to infinity''
(Corollary~\ref{cor:entropy-infinity}). These two results are proven
and used as tools to study interval maps in \cite{BR}, but they are 
interesting by themselves, that is why we state them here.

\begin{prop}\label{prop:entropy-infinity}
Let $G$ be a connected oriented graph which is not strongly positive
recurrent and $W$ a finite set of vertices. Then for all integers
$n$ there exists a connected subgraph $G_n\subset G$ such that $h(G_n)=h(G)$
and for all $w\in W$, for all $0\leq k<n$,
$f_{ww}^{G_n}(k)=0$.
\end{prop}

\begin{cor}\label{cor:entropy-infinity}
Let $G$ be a connected oriented graph which is not strongly positive
recurrent. Then there exists a sequence of ergodic Markov measures
$(\mu_n)_{n\geq 0}$ such that 
$\lim_{n\to+\infty}h_{\mu_n}(\Gamma_G,\sigma)=h(G)$ and for all finite
subsets of vertices $W$, $\disp\lim_{n\to+\infty}
\mu_n\left(\{(u_n)_{n\in\IZ}\in \Gamma_G\mid u_0\in W\}\right)=0$.
\end{cor}

\subsection{Local entropy and maximal measures}\label{subsec:hloc-mesure-max}
For a compact system, the {\em local entropy} is defined according to a 
distance but does not depend on it. One may wish to extend this 
definition to non compact metric spaces although the notion obtained in 
this way is not canonical.

\begin{defi}
Let $X$ be a metric space, $d$ its distance and let $T\colon X\to X$ be a
continuous map.

The {\em Bowen ball} of centre $x$, of radius $r$ and of order $n$ is defined
as 
$$B_n(x,r)=\{y\in X\mid d(T^ix,T^iy)<r , 0\leq i< n\}.$$

$E$ is a $(\delta,n)$-separated set if 
$$\forall y,y'\in E, y\not=y', \exists 0\leq k<n,\ d(T^k y, T^k y')\geq 
\delta.$$
The maximal cardinality of a $(\delta,n)$-separated set contained in $Y$ is
denoted by $s_n(\delta,Y)$.
 
The local entropy of $(X,T)$ is defined as 
$\disp h_{loc}(X)=\lim_{\eps\to 0}h_{loc}(X,\eps)$, where
$$
h_{loc}(X,\eps)=\lim_{\delta\to 0} \limsup_{n\to+\infty} 
\frac{1}{n}\sup_{x\in X}\log s_n(\delta,B_n(x,\eps)).
$$
\end{defi}

If the space $X$ is not compact, these notions depend on the distance.
When $X=\Gamma_G$, we use the distance
$d$ introduced in Section~\ref{subsec:markov-chain}. 
The local entropy of $\Gamma_G$ does not depend on the identification of
the vertices with $\IN$.

\begin{prop}
Let $\Gamma_G$ be the topological Markov chain on $G$ and
$\overline{\Gamma}_G$ its compactification as defined in 
Section~\ref{subsec:markov-chain}. Then $h_{loc}(\Gamma_G)=
h_{loc}(\overline{\Gamma}_G)$.
\end{prop}

\begin{demo}
Let $\bar u=(u_n)_{n\in\IZ} \in \overline{\Gamma}_G$, $\eps>0$ and $k\geq 1$.
By continuity there exists 
$\eta>0$ such that, if $\bar v\in \overline{\Gamma}_G$ and 
$d(\bar u,\bar v)<\eta$ then $d(\sigma^i(\bar u),\sigma^i(\bar v))<\eps$
for all $0\leq i<k$. By definition of $\overline{\Gamma}_G$ there is
$\bar v\in \Gamma_G$ such that $d(\bar u,\bar v)<\eta$, thus
$\bar u\in B_k(\bar v,\eps)$, which implies that 
$B_k(\bar u,\eps)\subset B_k(\bar v,2\eps)$. Consequently
$h_{loc}(\overline{\Gamma}_G,\eps)\leq h_{loc}(\Gamma,2\eps)$, and
$h_{loc}(\overline{\Gamma}_G)\leq h_{loc}(\Gamma_G)$. The reverse inequality
is obvious.
\end{demo}

We are going to prove that, if $h_{loc}(\Gamma_G)<h(G)$, then
$G$ is strongly positive recurrent. First we introduce some notations.

\medskip
Let $G$ be an oriented graph. If $V$ is a subset of vertices, $H$ a
subgraph of $G$
and $\bar{u}=(u_n)_{n\in\IZ}\in \Gamma_G$, define
$$
\CC^{H}(\bar u,V)=\{(v_n)_{n\in\IZ}\in \Gamma_H\mid
\forall n\in\IZ, u_n\in V\Rightarrow (v_n=u_n), 
u_n\not\in V\Rightarrow v_n\not\in V\}.
$$
If $S\subset \Gamma_G$ and $p,q\in\IZ\cup\{-\infty,+\infty\}$, define
$$
[S]_{p}^{q}=\{(v_n)_{n\in\IZ}\in \Gamma_G\mid
\exists (u_n)_{n\in\IZ}\in S,\forall p\leq n\leq q, u_n=v_n\}.
$$

\begin{lem}\label{lem:C-Bn}
Let $G$ be an oriented graph on the set of vertices $\IN$.
\begin{enumerate}
\item
If $V\supset \{0,\cdots,p+2\}$ then 
for all $\bar u\in\Gamma_G$ and all $n\geq 1$,
$\CC^G(\bar u,V)\subset B_n(\bar u,2^{-p}).$
\item
If $\bar u=(u_n)_{n\in\IZ}$ and $\bar v=(v_n)_{n\in\IZ}$
are two paths in $G$ such that $(u_0,\cdots,u_{n-1})
\not=(v_0,\cdots,v_{n-1})$ and $u_i,v_i\in\{0,\cdots,q-1\}$ for
$0\leq i\leq n-1$ then $(\bar u,\bar v)$ is $(2^{-q},n)$-separated.
\end{enumerate}
\end{lem}

\begin{demo}
\noindent (i)
Let $\bar u=(u_n)_{n\in\IZ}\in\Gamma_G$.
If $\bar v=(v_n)_{n\in\IZ}\in 
\CC^G(\bar u,V)$, then $D(u_j,v_j)\leq 2^{-(p+2)}$ for all $j\in\IZ$.
Consequently for all $0\leq i<n$
$$
d(\sigma^i(\bar u),\sigma^i(\bar v))
=\sum_{k\in\IZ}\frac{D(u_{i+k},v_{i+k})}{2^{|k|}}
\leq \sum_{k\in\IZ}\frac{2^{-(p+2)}}{2^{|k|}}
\leq  3\cdot 2^{-(p+2)}<  2^{-p}.
$$

\medskip
\noindent (ii)
Let $0\leq i\leq n-1$ such that $u_i\not= v_i$. By hypothesis, 
$u_i,v_i\leq q-1$. Suppose that $u_i<v_i$. Then
$\disp
d(\sigma^i (\bar u),\sigma^i (\bar v))\geq D(u_i,v_i)=
2^{-u_i}(1-2^{-(v_i-u_i)})\geq 2^{-q}$.
\end{demo}

\begin{theo}\label{theo:hloc-maximal-measure}
Let $G$ be a connected oriented graph of finite entropy on the set of
vertices $\IN$.
If $h_{loc}(\Gamma_G)<h(G)$, then the graph $G$ is strongly positive recurrent
and the Markov chain $(\Gamma_G,\sigma)$ 
admits a maximal measure.
\end{theo}

\begin{demo}
Fix $C$ and $\eps>0$ such that $h_{loc}(\Gamma_G,\eps)<C<h(G)$.
Let $p$ be an integer such that $2^{-(p-1)}<\eps$. Let $G'$ be
a finite subgraph $G'$ such that
$h(G')>C$ and let $V$ be a finite subset of vertices
such that $V$ is connected and contains the vertices of $G'$ and the vertices
$\{0,\cdots,p\}$. Define $W=\overline{V}$,
$V_q=\{n\leq q\}$ and $W_q=V_q\setminus V=W\cap V_q$ 
for all $q\geq 1$.

Our aim is to bound $t_{uu'}^{W}(n)=t_{uu'}^{\overline {V}}(n)$.
Choose $u,u'\in V$ and let $(w_0,\cdots,w_{n_0})$ be a path between $u'$
and $u$ with $w_i\in V$ for $0\leq i\leq n_0$. Fix $n\geq 1$. 
One has $\disp t_{uu'}^W(n)=\lim_{q\to+\infty} t_{uu'}^{W_q}(n)$.

Fix $\delta_0>0$ such that 
$$
\forall \delta\leq\delta_0,\limsup_{n\to+\infty}\frac{1}{n}
\sup_{\bar v\in\Gamma_G}\log s_n(\delta,B_n(\bar v,\eps))<C.
$$
Take $q\geq 1$ arbitrarily large and $\delta\leq \min\{\delta_0,
2^{-(q+1)}\}$. Choose $N$ such that 
\begin{equation}
\forall n\geq N, \forall \bar v\in\Gamma_G, \frac{1}{n}\log 
s_n(\delta,B_n(\bar v,\eps))<C.
\label{eq:N} \end{equation}

If $t_{uu'}^{W_q}(n)\not=0$, choose a path $(v_0,\cdots,v_n)$ 
such that $v_0=u, v_n=u'$ and $v_i\in W_q$ for $0<i<q$.
Define $\bar v^{(n)}=(v^{(n)}_i)_{i\in\IZ}$ as the periodic path of period
$n+n_0$ satisfying
$v^{(n)}_i=v_i$ for $0\leq i\leq n$  and
$v^{(n)}_{n+i}=w_i$ for $0\leq i\leq n_0$.

Define the set $E_q(n,k)$ as follows (see Figure~\ref{fig:Eq}):
$$ 
E_q(n,k)=\left[\CC^{V_q}(\bar v^{(n)},V)\right]_0^{k(n+n_0)}\cap
\left[\bar v^{(n)}\right]_{-\infty}^{0}\cap
\left[\bar v^{(n)}\right]_{k(n+n_0)}^{+\infty}.
$$

\begin{figure}[ht]
\centerline{\includegraphics{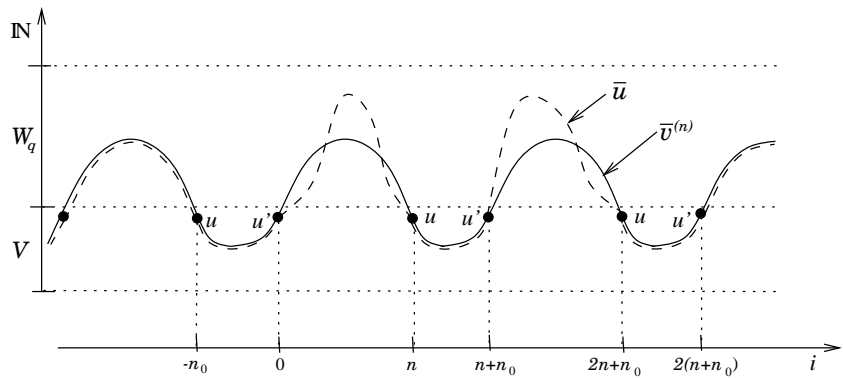}}
\caption{The set $E_q(n,k)$ ($k=2$ on the picture): 
$\bar{v}^{(n)}$ (in solid) is a periodic path, 
$\bar u$ (in dashes) is a element of  $E_q(n,k)$. Between the indices
$0$ and $k(n+n_0)$, $\bar v^{(n)}$ and $\bar u$ coincide when $v^{(n)}_i$
is in $V$ and $\bar v^{(n)}$ and $\bar u$ are in $W_q$ at the same time.
Before $0$ or after $k(n+n_0)$, the two paths coincide\label{fig:Eq}.}
\end{figure}

The paths in $E_q(n,1)$ are exactly the paths 
counted by $t_{uu'}^{W_q}(n)$ which are extended outside the indices
$\{0,\cdots,n\}$ like the path $\bar v^{(n)}$, 
thus $\# E_q(n,1)=t_{uu'}^{W_q}(n)$.
Similarly, $\# E_q(n,k)=\left(t_{uu'}^{W_q}(n)\right)^k$.

By definition, $E_q(n,k)\subset \CC^G(\bar v^{(n)},V)$ and 
$\{0,\cdots,p\}\subset V$ thus  
$E_q(n,k)\subset B_{k(n+n_0)}(\bar v^{(n)},\eps)$ by Lemma~\ref{lem:C-Bn}(i).
Moreover, if $(w_i)_{i\in\IZ}$ and $(w_i')_{i\in\IZ}$ are two distinct
elements of $E_q(n,k)$, there exists $0\leq i <k(n+n_0)$ such that 
$w_i\not=w_i'$ and $w_i,w_i'\leq q$, thus $E_q(n,k)$ 
is a $(\delta,k(n+n_0))$-separated set by Lemma~\ref{lem:C-Bn}(ii).
Choose $k$ such that $k(n+n_0)\geq N$. Then by Equation~(\ref{eq:N})
$$
\# E_q(n,k)\leq s_{k(n+n_0)}(\delta,B_{k(n+n_0)}(\bar v^{(n)},\eps))
< e^{k(n+n_0)C}.
$$
As $\# E_q(n,k)=\left(t_{uu'}^{W_q}(n)\right)^k$, one gets
$t_{uu'}^{W_q}(n)<e^{(n+n_0)C}$. 
This is true for all $q\geq 1$, thus 
$$
t_{uu'}^W(n)=\lim_{q\to +\infty}t_{uu'}^{W_q}(n)\leq e^{(n+n_0)C}
$$
and
$$
\tau_{uu'}^W=\tau_{uv'}^{\overline{V}}\leq C<h(V).
$$
Theorem~\ref{theo:GZ-et-reciproque} concludes the proof.
\end{demo}

\begin{rem}
Define the entropy at infinity as
$h_{\infty}(G)=\lim_{n\to+\infty}h(G\setminus G_n)$
where $(G_n)_{n\geq 0}$ is a sequence of finite graphs such that
$\bigcup_n G_n=G$. The local entropy satisfies 
$h_{loc}(\Gamma_G)\geq h_{\infty}(G)$ but in general these two quantities
are not equal and the condition $h_{\infty}(G)<h(G)$ does
not imply that $G$ is strongly positive recurrent. This is illustrated by
Example~\ref{ex:R=L} (see Figure~\ref{fig:G-G'}).
\end{rem}


\noindent
Sylvie {\sc Ruette} -- Institut de Math\'ematiques de Luminy -- 
CNRS UPR 9016 -- 163 avenue de Luminy, case 907 -- 13288 Marseille cedex 9 
-- France\\
{\it E-mail :} {\tt  ruette@iml.univ-mrs.fr}

\end{document}